\input amstex
\documentstyle{amsppt}
%
\catcode`@=11
\redefine\output@{%
  \def\break{\penalty-\@M}\let\par\endgraf
  \ifodd\pageno\global\hoffset=105pt\else\global\hoffset=8pt\fi  
  \shipout\vbox{%
    \ifplain@
      \let\makeheadline\relax \let\makefootline\relax
    \else
      \iffirstpage@ \global\firstpage@false
        \let\rightheadline\frheadline
        \let\leftheadline\flheadline
      \else
        \ifrunheads@ 
        \else \let\makeheadline\relax
        \fi
      \fi
    \fi
    \makeheadline \pagebody \makefootline}%
  \advancepageno \ifnum\outputpenalty>-\@MM\else\dosupereject\fi
}
\def\Beta{\mathchar"0\hexnumber@\rmfam 42}
\catcode`\@=\active
\nopagenumbers
\def\negskp{\hskip -2pt}

\def\const{\operatorname{const}}
\def\blue#1{#1}

\catcode`#=11\def\diez{#}\catcode`#=6
\catcode`_=11\def\podcherkivanie{_}\catcode`_=8
\def\mycite#1{\cite{\blue{#1}}\immediate\special{ps:
     ShrHPSdict begin /ShrBORDERthickness 0 def}}

\def\mytag#1{%
    \tag#1}
\def\mythetag#1{\thetag{\blue{#1}}\immediate\special{ps:
     ShrHPSdict begin /ShrBORDERthickness 0 def}}
\def\myrefno#1{\no#1}
\def\myhref#1#2{\blue{#2}\immediate\special{ps:
     ShrHPSdict begin /ShrBORDERthickness 0 def}}
\def\myEarXivlink{\myhref{http://arXiv.org}{http:/\negskp/arXiv.org}}

\def\mytheorem#1{\csname proclaim\endcsname{Theorem #1}}
\def\mytheoremwithtitle#1#2{\csname proclaim\endcsname{Theorem #1#2}}
\def\mythetheorem#1{\blue{#1}\immediate\special{ps:
     ShrHPSdict begin /ShrBORDERthickness 0 def}}
\def\mylemma#1{\csname proclaim\endcsname{Lemma #1}}
\def\mylemmawithtitle#1#2{\csname proclaim\endcsname{Lemma #1#2}}
\def\mythelemma#1{\blue{#1}\immediate\special{ps:
     ShrHPSdict begin /ShrBORDERthickness 0 def}}
\def\mycorollary#1{\csname proclaim\endcsname{Corollary #1}}

\def\myconjecture#1{\csname proclaim\endcsname{Conjecture #1}}
\def\myconjecturewithtitle#1#2{\csname proclaim\endcsname{Conjecture #1#2}}
\def\mytheconjecture#1{\blue{#1}\immediate\special{ps:
     ShrHPSdict begin /ShrBORDERthickness 0 def}}

\pagewidth{360pt}
\pageheight{606pt}
\topmatter
\title
A note on Khabibullin's conjecture for integral inequalities.
\endtitle
\author
R.~A.~Sharipov
\endauthor
\address 5 Rabochaya street, 450003 Ufa, Russia\newline
\vphantom{a}\kern 12pt Cell Phone: +7(917)476 93 48
\endaddress
\email \vtop to 30pt{\hsize=280pt\noindent
\myhref{mailto:r-sharipov\@mail.ru}
{r-sharipov\@mail.ru}\newline
\myhref{mailto:R\podcherkivanie Sharipov\@ic.bashedu.ru}
{R\_\hskip 1pt Sharipov\@ic.bashedu.ru}\vss}
\endemail
\urladdr
\vtop to 20pt{\hsize=280pt\noindent
\myhref{http://ruslan-sharipov.ucoz.com/}
{http:/\negskp/ruslan-sharipov.ucoz.com}\newline
\myhref{http://www.freetextbooks.narod.ru/}
{http:/\negskp/www.freetextbooks.narod.ru}\newline
\myhref{http://sovlit2.narod.ru/}
{http:/\negskp/sovlit2.narod.ru}\vss}
\endurladdr
\abstract
    An integral transformation relating two inequalities in Khabibullin's 
conjecture is found. Another proof of this conjecture for some special values 
of its numeric parameters is suggested. 
\endabstract
\endtopmatter
\TagsOnRight
\document

\head
1. Introduction.
\endhead
\myconjecturewithtitle{1.1}{ (Khabibullin)} Let $\alpha>1/2$ and let $q=q(t)$
be a positive continuous function on the half-line $[\kern 1pt 0,+\infty)$. 
Then the inequality 
$$
\hskip -2em
\int\limits^{\,1}_0\left(\,\,\int\limits^{\,1}_x(1-y)^{n-1}\,\frac{dy}{y}\right)
q(tx)\,dx\leqslant t^{\alpha-1}
\mytag{1.1}
$$
fulfilled for all \kern 2pt $0\leqslant t<+\infty$ implies the inequality
$$
\hskip -2em
\int\limits^{+\infty}_0 q(t)\,\ln\Bigl(1+\frac{1}{t^{\,2\kern 0.2pt\alpha}}\Bigr)\,dt
\leqslant\pi\,\alpha\prod^{n-1}_{k=1}\Bigl(1+\frac{\alpha}{k}\Bigr).
\mytag{1.2}
$$
\endproclaim
     The conjecture~\mytheconjecture{1.1} arose in \mycite{1} (see also \mycite{2}), 
though in some different form. The statement of this conjecture in the above form is 
given in \mycite{3}.\par
     Note that the conjecture~\mytheconjecture{1.1} is formulated for $\alpha>1/2$.
Actually it could be formulated for all positive $\alpha$ in the following way. 
\myconjecturewithtitle{1.2}{ (Khabibullin)} Let $\alpha>0$ and let $q=q(t)$
be a positive conti\-nuous function on the half-line $[\kern 1pt 0,+\infty)$. 
Then the inequality 
$$
\hskip -2em
\int\limits^{\,1}_0\left(\,\,\int\limits^{\,1}_x(1-y)^{n-1}\,\frac{dy}{y}\right)
q(tx)\,dx\leqslant t^{\alpha-1}
\mytag{1.3}
$$
fulfilled for all \kern 2pt $0\leqslant t<+\infty$ implies the inequality
$$
\pagebreak
\hskip -2em
\int\limits^{+\infty}_0 q(t)\,\ln\Bigl(1+\frac{1}{t^{\,2\kern 0.2pt\alpha}}\Bigr)\,dt
\leqslant\pi\,\alpha\prod^{n-1}_{k=1}\Bigl(1+\frac{\alpha}{k}\Bigr).
\mytag{1.4}
$$
\endproclaim
     The matter is that in \mycite{4} the conjecture~\mytheconjecture{1.2} is already 
proved to be valid for $0<\alpha\leqslant 1/2$. Despite this fact, in the present paper 
we shall use the conjecture~\mytheconjecture{1.1} in its more general 
form~\mytheconjecture{1.2}. The main goal of this paper is to give another treatment of 
the inequalities \mythetag{1.1}, \mythetag{1.2}, \mythetag{1.3}, \mythetag{1.4} and to 
give another proof of the conjecture~\mytheconjecture{1.2} for the case $0<\alpha
\leqslant 1/2$. Some other particular values of the parameters $n$ and $\alpha$ are
also considered. 
\head
2. Relation to Euler's Beta function. 
\endhead
     The product in the right hand side of the inequalities \mythetag{1.2} and
\mythetag{1.4} is related to Euler's Beta function. We have the following relationship:
$$
\hskip -2em
\alpha\,\prod^{n-1}_{k=1}\Bigl(1+\frac{\alpha}{k}\Bigr)=\frac{1}{\Beta(\alpha,n)}.
\mytag{2.1}
$$
Let's recall that Euler's Beta function $\Beta(\alpha,\beta)$ is defined through the 
integral
$$
\hskip -2em
\Beta(\alpha,\beta)=\int\limits^{\,1}_0x^{\alpha-1}\,(1-x)^{\beta-1}\,dx.
\mytag{2.2}
$$
The Beta function \mythetag{2.2} is expressed through Euler's Gamma function
$$
\hskip -2em
\Beta(\alpha,\beta)=\frac{\Gamma(\alpha)\,\Gamma(\beta)}{\Gamma(\alpha+\beta)}.
\mytag{2.3}
$$
The value of the Gamma function $\Gamma(\beta)$ in \mythetag{2.3} is expressed through the 
factorial for integer values of the argument $\beta=n$:
$$
\hskip -2em
\Gamma(n)=(n-1)!\text{ \ for all \ }n=1,\,2,\,3\,\ldots\,.
\mytag{2.4}
$$
Applying \mythetag{2.4} to \mythetag{2.3}, we recall the following property of the
Gamma function:
$$
\hskip -2em
\Gamma(\alpha)=(\alpha-1)\,\Gamma(\alpha-1).
\mytag{2.5}
$$
Now from \mythetag{2.3}, \mythetag{2.4}, and \mythetag{2.5} we derive the formula
$$
\hskip -2em
\Beta(\alpha,n)=\frac{(n-1)!}{(\alpha+n-1)\,(\alpha+n-2)\cdot\ldots\cdot(\alpha+1)\,\alpha}.
\mytag{2.6}
$$
The formula \mythetag{2.6} can be transformed as follows:
$$
\hskip -2em
\Beta(\alpha,n)=\frac{1}{\alpha}\prod^{n-1}_{k=1}\frac{k}{k+\alpha}
=\frac{1}{\alpha}\prod^{n-1}_{k=1}\frac{1}{\lower 2pt\hbox{$1
+\dfrac{\lower 2pt\hbox{$\alpha$}}{k}$}}.
\mytag{2.7}
$$
Now it is easy to see that the formula \mythetag{2.7} is equivalent to the relationship 
\mythetag{2.1}. Thus, the formula \mythetag{2.1} is proved.\par
     Note that the formula \mythetag{2.3} as well as the formulas equivalent to 
\mythetag{2.4} and \mythetag{2.5} can be found in 
Chapter~\uppercase\expandafter{\romannumeral 17}  of the book \mycite{5}. The proofs of
these formulas are also available there. 
\head
3. The study of the kernel.
\endhead
     Let's denote through $A_n(x)$ the following integral:
$$
\hskip -2em
A_n(x)=\int\limits^{\,1}_{\!x}(1-y)^n\,\frac{dy}{y}.
\mytag{3.1}
$$
The integral \mythetag{3.1} is used as a kernel in the integral inequality 
\mythetag{1.3}. In terms of our notation $A_n(x)$ the inequality \mythetag{1.3}
looks like
$$
\hskip -2em
\int\limits^{\,1}_0 A_{n-1}(x)\,q(tx)\,dx\leqslant t^{\alpha-1}.
\mytag{3.2}
$$\par
     The integral \mythetag{3.1} can be calculated explicitly. We use the relationship
$$
\hskip -2em
A_n(x)=A_{n-1}(x)-\frac{(1-x)^n}{n}\text{\ \ for \ }n\geqslant 1.
\mytag{3.3}
$$
The recurrent relationship is derived immediately from \mythetag{3.1}. Indeed, we have
$$
\gather
A_n(x)=\int\limits^{\,1}_{\!x}(1-y)\,(1-y)^{n-1}\,\frac{dy}{y}
=\int\limits^{\,1}_{\!x}(1-y)^{n-1}\,\frac{dy}{y}\,-\\
-\int\limits^{\,1}_{\!x}(1-y)^{n-1}\,dy=A_{n-1}(x)+
\frac{(1-y)^n}{n}\,{\vrule height 16pt depth 7pt}^{\,1}_{\,x}
=A_{n-1}(x)-\frac{(1-x)^n}{n}.
\endgather
$$
Applying the relationship \mythetag{3.3} recursively, we derive
$$
\hskip -2em
A_n(x)=A_0(x)-\sum^n_{m=1}\frac{(1-x)^m}{m}.
\mytag{3.4}
$$
Note that the term $A_0(x)$ in \mythetag{3.3} is calculated explicitly. Indeed, 
we have
$$
\hskip -2em
A_0(x)=\int\limits^{\,1}_{\!x}\frac{dy}{y}=\ln y\,{\vrule height 13pt depth 
7pt}^{\,1}_{\,x}=-\ln x.
\mytag{3.5}
$$
\mylemma{3.1} The kernel $A_n(x)$ is given by the explicit formula
$$
\hskip -2em
A_n(x)=-\ln x-\sum^n_{m=1}\frac{(1-x)^m}{m}.
\mytag{3.6}
$$
\endproclaim
    The proof of the lemma~\mythelemma{3.1} is immediate from \mythetag{3.4} and 
\mythetag{3.5}.
\mylemma{3.2} The kernel $A_n(x)$ is a continuous function on the interval
$(0,1]$ with a logarithmic singularity at the the point $x=0$. It vanishes 
at the point $x=1$.
\endproclaim
    The proof of the lemma~\mythelemma{3.2} is immediate from the formula
\pagebreak \mythetag{3.6}.\par
    Note that the logarithmic function has the following Tailor expansion:
$$
\hskip -2em
-\ln x=-\ln(1-(1-x))=\sum^{\infty}_{m=1}\frac{(1-x)^m}{m}.
\mytag{3.7}
$$
Combining \mythetag{3.7} with the formula \mythetag{3.6}, we get
$$
\hskip -2em
A_n(x)=\kern -4pt\sum^\infty_{m=n+1}\kern -4pt\frac{(1-x)^m}{m}.
\mytag{3.8}
$$
\mylemma{3.3}\parshape 3 0pt 360pt 0pt 360pt 180pt 180pt The kernel $A_n(x)$ is a 
positive decreasing function on the interval $(0,1)$ vanishing at the point $x=1$. 
The kernel function $A_n(x)$ and its derivatives $A'_n(x),\,A''_n(x),\,\ldots,
\,A^{(n)}_n(x)$ up to the $n$-th order do vanish at the point $x=1$.
\endproclaim
\demo{Proof}\parshape 1 180pt 180pt  The power series \mythetag{3.8} converges 
on the interval $(0,1]$. 
\vadjust{\vskip 5pt\hbox to 0pt{\kern 0pt \includegraphics{khab01.eps}\hss}
\vskip -5pt}Each term of this series is a positive decreasing function on the interval 
$(0,1]$. Hence the kernel $A_n(x)$ is a positive decreasing function on the interval 
$(0,1]$.\par
\parshape 1 180pt 180pt The equality $A_n(1)=0$ and the equalities $A^{(k)}_n(1)=0$
for $k=1,\,\ldots,\,n$ are derived from the formula \mythetag{3.8}. 
\qed\enddemo
\parshape 3 180pt 180pt 180pt 180pt 0pt 360pt Relying on the lemmas \mythelemma{3.2} 
and \mythelemma{3.3}, one can plot the graph of the kernel function $A_n(x)$. 
This graph is shown above on Fig\.~3.1.\par
     Since all terms of the power series \mythetag{3.8} are positive on the interval
$(0,1)$, we derive the following inequalities for the kernel functions on this interval:
$$
\hskip -2em
A_0(x)>A_1(x)>A_2(x)>\ldots>A_n(x)>\ldots\,.
\mytag{3.9}
$$
These kernel functions \mythetag{3.9} vanish simultaneously at the point $x=1$.
\head
4. The basic example.
\endhead
     Let's consider the function $q(x)=C\,x^{\alpha-1}$, where $\alpha>1/2$ and
$C=\const>0$. If $\alpha<1$, this function is not continuous at $x=0$. Nevertheless, 
the integral 
$$
\hskip -2em
\int\limits^{\,1}_0\! A_{n-1}(x)\,q(tx)\,dx=C\,t^{\alpha-1}\!\!\int\limits^{\,1}_0 
\!A_{n-1}(x)\,x^{\alpha-1}\,dx
\mytag{4.1}
$$
converges at the point $x=0$ for any $\alpha>0$. Relying on \mythetag{4.1}, we denote 
$$
\hskip -2em
I_{n\,\alpha}=\int\limits^{\,1}_0\!A_{n-1}(x)\,x^{\alpha-1}\,dx.
\mytag{4.2}
$$
The constant $I_{n\,\alpha}$ can be calculated explicitly. Substituting \mythetag{3.6}
into \mythetag{4.2}, we get
$$
\hskip -2em
I_{n\,\alpha}=-\int\limits^{\,1}_0\!\ln(x)\,x^{\alpha-1}\,dx
-\sum^{n-1}_{m=1}\int\limits^{\,1}_0\frac{x^{\alpha-1}\,(1-x)^m}{m}\,dx.
\mytag{4.3}
$$
For $\alpha>0$ the first integral from \mythetag{4.3} is calculated as follows:
$$
\hskip -2em
\int\limits^{\,1}_0\!\ln(x)\,x^{\alpha-1}\,dx=\ln(x)\,\frac{x^\alpha}{\alpha}
\,{\vrule height 16pt depth 7pt}^{\kern 1.5pt 1}_{\kern 2pt 0}
-\int\limits^{\,1}_0\!\frac{x^{\alpha-1}}{\alpha}\,dx=-\frac{x^\alpha}{\alpha^2}
\,{\vrule height 16pt depth 7pt}^{\kern 1.5pt 1}_{\kern 2pt 0}=-\frac{1}{\alpha^2}.
\mytag{4.4}
$$
Due to \mythetag{2.2} the second integral from \mythetag{4.3} is expressed through 
the Beta function:
$$
\hskip -2em
\int\limits^{\,1}_0 x^{\alpha-1}\,(1-x)^m\,dx=B(\alpha,m+1).
\mytag{4.5}
$$
Combining \mythetag{4.4} and \mythetag{4.5}, we derive
$$
\hskip -2em
I_{n\,\alpha}=\frac{1}{\alpha^2}-\sum^{n-1}_{m=1}\frac{B(\alpha,m+1)}{m}.
\mytag{4.6}
$$\par
      Note that there is another way for calculating the constant $I_{n\,\alpha}$. 
From the formula \mythetag{4.2} we immediately derive the following expression for
$I_{n\,\alpha}$:
$$
\hskip -2em
I_{n\,\alpha}=A_{n-1}(x)\,\frac{x^\alpha}{\alpha}
\,{\vrule height 16pt depth 7pt}^{\kern 1.5pt 1}_{\kern 2pt 0}
-\int\limits^{\,1}_0\frac{A'_{n-1}(x)\,x^\alpha}{\alpha}\,dx.
\mytag{4.7}
$$
Since $\alpha>0$ and $A_{n-1}(1)=0$ for each integer $n\geqslant 1$, the formula
\mythetag{4.7} yields 
$$
\hskip -2em
I_{n\,\alpha}=-\int\limits^{\,1}_0\frac{A'_{n-1}(x)\,x^\alpha}{\alpha}\,dx=
\int\limits^{\,1}_0\frac{x^{\alpha-1}\,(1-x)^{n-1}}{\alpha}\,dx
=\frac{B(\alpha,n)}{\alpha}.
\mytag{4.8}
$$
Comparing \mythetag{4.8} and \mythetag{4.6}, for each integer $n\geqslant 1$ we 
derive the identity
$$
\hskip -2em
\frac{B(\alpha,n)}{\alpha}=\frac{1}{\alpha^2}-\sum^{n-1}_{m=1}
\frac{B(\alpha,m+1)}{m}.
\mytag{4.9}
$$\par
     For $n=1$ the identity \mythetag{4.9} simplifies. In this case it is written 
as follows:
$$
\hskip -2em
B(\alpha,1)=\frac{1}{\alpha}.
\mytag{4.10}
$$
The equality \mythetag{4.10} is valid for $\alpha>0$. It is easily derived immediately 
from the formula \mythetag{2.2} defining Euler's Beta function.\par
     For $\alpha>0$ and for each integer $n\geqslant 2$, applying the formula \mythetag{2.7} 
to \mythetag{4.9}, we transform the identity \mythetag{4.9} to the following one:
$$
\hskip -2em
\frac{1}{\alpha}\prod^{n-1}_{k=1}\frac{k}{k+\alpha}=\frac{1}{\alpha}-\sum^{n-1}_{m=1}
\!\left(\,\prod^m_{k=1}\frac{k}{m\,(k+\alpha)}\right).
\mytag{4.11}
$$
For $n=2$ the identity \mythetag{4.11} is verified by means of direct calculations. Then
for each integer $n>2$ it is proved by induction on $n$.\par
      Now, returning back to the function $q(x)=C\,x^{\alpha-1}$, we choose $C=1/I_{n\,\alpha}$.
Under this special choice of the constant $C$ the function 
$$
\hskip -2em
q(x)=\frac{x^{\alpha-1}}{I_{n\,\alpha}}=\frac{\alpha\,x^{\alpha-1}}{B(\alpha,n)}
\mytag{4.12}
$$
turns the inequality \mythetag{3.2} into the equality
$$
\hskip -2em
\int\limits^{\,1}_0 A_{n-1}(x)\,q(tx)\,dx=t^{\alpha-1}.
\mytag{4.13}
$$
Let's substitute the function \mythetag{4.12} into the integral in the left hand side
of \mythetag{1.4}:
$$
\hskip -2em
\int\limits^{+\infty}_0 q(t)\,\ln\Bigl(1+\frac{1}{t^{\,2\kern 0.2pt\alpha}}\Bigr)\,dt
=\frac{\alpha}{B(\alpha,n)}\int\limits^{+\infty}_0\!t^{\alpha-1}\,\ln\Bigl(1+\frac{1}
{t^{\,2\kern 0.2pt\alpha}}\Bigr)\,dt.
\mytag{4.14}
$$
The integral in the right hand side of \mythetag{4.14} is calculated explicitly as an
indefinite integral. Indeed, by differentiation one can verify that
$$
\hskip -2em
\int\limits^{\,x}\!t^{\alpha-1}\,\ln\Bigl(1+\frac{1}
{t^{\,2\kern 0.2pt\alpha}}\Bigr)\,dt=\frac{t^\alpha}{\alpha}\,\ln\Bigl(1+\frac{1}
{t^{\,2\kern 0.2pt\alpha}}\Bigr)-\frac{2}{\alpha}\,\arctan(t^{-\alpha}).
\mytag{4.15}
$$
Assuming that $\alpha>0$ and applying the formula \mythetag{4.15}, we calculate the 
integral
$$
\hskip -2em
\int\limits^{+\infty}_0\!t^{\alpha-1}\,\ln\Bigl(1+\frac{1}
{t^{\,2\kern 0.2pt\alpha}}\Bigr)\,dt=\frac{\pi}{\alpha}.
\mytag{4.16}
$$
Now from \mythetag{4.14} and \mythetag{4.16} for the function \mythetag{4.12}
we derive the equality 
$$
\hskip -2em
\int\limits^{+\infty}_0 q(t)\,\ln\Bigl(1+\frac{1}{t^{\,2\kern 0.2pt\alpha}}\Bigr)\,dt
=\frac{\pi}{B(\alpha,n)}.
\mytag{4.17}
$$
Combining \mythetag{4.17} with \mythetag{2.1}, we find that the inequality \mythetag{1.4}
turns to the equality for the function \mythetag{4.12}. Indeed, we have
$$
\pagebreak
\hskip -2em
\int\limits^{+\infty}_0 q(t)\,\ln\Bigl(1+\frac{1}{t^{\,2\kern 0.2pt\alpha}}\Bigr)\,dt
=\pi\,\alpha\prod^{n-1}_{k=1}\Bigl(1+\frac{\alpha}{k}\Bigr).
\mytag{4.18}
$$
\mylemma{4.1} The special choice of the function $q(x)$ given by the formula
\mythetag{4.12} turn both inequalities \mythetag{1.3} and \mythetag{1.4} into
equalities. 
\endproclaim
     The lemma~\mythelemma{4.1} is immediate from \mythetag{4.13} and \mythetag{4.18}.
This lemma could be a motivation for the conjecture~\mytheconjecture{1.2}.
\head
5. Some integral relationships.
\endhead
     Let $\varphi(t)$ be a smooth function defined on the interval $(0,+\infty)$ 
and such that it satisfies the following asymptotic condition for some $\varepsilon>0$:
$$
\hskip -2em
\varphi(t)=O(t^{-\varepsilon})\text{\ \ as \ }t\to +\infty.
\mytag{5.1}
$$
In addition to \mythetag{5.1}, assume that its derivatives satisfy the conditions
$$
\hskip -2em
\varphi^{(s)}(t)=\frac{d^s\varphi}{dt^s}=O(t^{-s-\varepsilon})\text{\ \ as \ }t\to +\infty
\text{\ \ for all \ }s=1,\,\ldots,\,n+2.
\mytag{5.2}
$$
By means of the function $\varphi(t)$ we define the function
$$
\hskip -2em
\varPhi_n(t)=-\frac{d}{dt}\left(\frac{(-t)^{n+1}}{n!}\,\,\varphi^{(n+1)}(t)\right)
\!\text{, \ where \ } \varphi^{(n+1)}(t)=\frac{d^{\,n+1}\varphi}{dt^{n+1}}.
\mytag{5.3}
$$
\mylemma{5.1} The function \mythetag{5.3} satisfies the integral relationship
$$
\hskip -2em
\int\limits^{+\infty}_y\!\!\varPhi_n(t)\,A_n(y/t)\,dt=\varphi(y)
\text{\ \ for \ }n\geqslant 0,
\mytag{5.4}
$$
provided the conditions \mythetag{5.1} and \mythetag{5.2} are fulfilled. 
\endproclaim
     The kernel function $A_n(x)$ has the logarithmic singularity at the point $x=0$ 
(see Lemma~\mythelemma{3.2}). Indeed, from \mythetag{3.6} we derive $A_n(x)\sim -\ln(x)$
as $x\to 0$. Then the function $A_n(y/t)$ has the logarithmic singularity 
at the infinity:
$$
\hskip -2em
A_n(y/t)\sim \ln t\text{\ \ as \ }t\to +\infty.
\mytag{5.5}
$$
On the other hand, from \mythetag{5.1} and \mythetag{5.2} for the function
\mythetag{5.3} we derive 
$$
\hskip -2em
\varPhi_n(t)=O(t^{-1-\varepsilon})\text{\ \ as \ }t\to +\infty.
\mytag{5.6}
$$
Combining \mythetag{5.5} with \mythetag{5.6}, we find that the integral
\mythetag{5.4} converges at infinity. Since the product $\varPhi_n(t)\,A_n(y/t)$
has no singularities at finite points $t\in[\,y,+\infty)$, the integral \mythetag{5.4}
converges in whole. 
\demo{Proof of the lemma~\mythelemma{5.1}} The proof is pure calculations. Upon 
substituting \mythetag{5.3} into the integral \mythetag{5.4} we can integrate by 
parts:
$$
\allowdisplaybreaks
\gather
\int\limits^{+\infty}_y \varPhi_n(t)\,A_n(y/t)\,dt=-\!\!\int\limits^{+\infty}_y 
\frac{d}{dt}\left(\frac{(-t)^{n+1}}{n!}\,\,\varphi^{(n+1)}(t)\right)\!A_n(y/t)
\,dt=\\
=-\frac{(-t)^{n+1}}{n!}\,\,\varphi^{(n+1)}(t)\,A_n(y/t)
\,{\vrule height 16pt depth 7pt}^{\,\infty}_{\ y}-\int\limits^{+\infty}_y 
\frac{(-t)^{n+1}}{n!}\,\,\varphi^{(n+1)}(t)\,A'_n(y/t)\,\frac{y\,dt}{t^2}.
\endgather
$$
From \mythetag{5.2} we derive $(-t)^{n+1}\,\varphi^{(n+1)}(t)=O(t^{-\varepsilon})$ as
$t\to +\infty$. Therefore the product $(-t)^{n+1}\,\varphi^{(n+1)}(t)$ suppresses the
logarithmic singularity of the function $A_n(y/t)$ at infinity. As for the lower
limit $t=y$, the function $A_n(y/t)$ vanishes at this point. Indeed, we have
$A_n(y/y)=A_n(1)=0$ (see Lemma~\mythelemma{3.2} or Lemma~\mythelemma{3.3}). As a result 
the integral \mythetag{5.4} reduces to the following one:
$$
\hskip -2em
\int\limits^{+\infty}_y\!\!\varPhi_n(t)\,A_n(y/t)\,dt=-\!\!\int\limits^{+\infty}_y 
\frac{(-t)^{n-1}}{n!}\,\,\varphi^{(n+1)}(t)\,A'_n(y/t)\,y\,dt.
\mytag{5.7}
$$
The first derivative $A'_n(y/t)$ in \mythetag{5.7} can be calculated explicitly.
Indeed, differentiating the formula \mythetag{3.6}, we derive the following expression 
for $A'_n(x)$:
$$
\hskip -2em
A'_n(x)=-\frac{1}{x}+\sum^n_{m=1}(1-x)^{m-1}.
\mytag{5.8}
$$
The sum in \mythetag{5.8} is the sum of a geometric progression. It is calculated 
explicitly:
$$
\hskip -2em
A'_n(x)=-\frac{1}{x}+\frac{(1-x)^n-1}{(1-x)-1}=\frac{-(1-x)^n}{x}.
\mytag{5.9}
$$
Substituting $x=y/t$ into the formula \mythetag{5.9}, we get
$$
\hskip -2em
A'_n(y/t)=\frac{-(1-y/t)^n}{y/t}=\frac{-(t-y)^n}{t^{n-1}\,y}.
\mytag{5.10}
$$
The next step is to substitute \mythetag{5.10} into \mythetag{5.7}. As a result
we derive
$$
\hskip -2em
\int\limits^{+\infty}_y\!\!\varPhi_n(t)\,A_n(y/t)\,dt=-\!\!\int\limits^{+\infty}_y 
\frac{(y-t)^n}{n!}\,\,\varphi^{(n+1)}(t)\,\,dt.
\mytag{5.11}
$$\par
     The integral in the right hand side of the formula can be calculated by means
of integrating by parts. Indeed, we easily derive the formula
$$
\int\limits^{+\infty}_y\!\!\varPhi_n(t)\,A_n(y/t)\,dt=-\frac{(y-t)^n}{n!}
\,\,\varphi^{(n)}(t)\,{\vrule height 16pt depth 7pt}^{\,\infty}_{\ y}
-\int\limits^{+\infty}_y\frac{(y-t)^{n-1}}{(n-1)!}\,\,\varphi^{(n)}(t)\,\,dt.
$$
Note that $(y-t)^n\,\,\varphi^{(n)}(t)=O(t^{-\varepsilon})$ as $t\to +\infty$ due
to \mythetag{5.2}. Moreover, $(y-t)^n=0$ if $n>0$ and $t=y$. Therefore the above
formula reduces to 
$$
\hskip -2em
\int\limits^{+\infty}_y\!\!\varPhi_n(t)\,A_n(y/t)\,dt=-\!\!\int\limits^{+\infty}_y
\frac{(y-t)^{n-1}}{(n-1)!}\,\,\varphi^{(n)}(t)\,\,dt.
\mytag{5.12}
$$
Comparing \mythetag{5.11} and \mythetag{5.12}, we see that the above calculations
let us to pass from $n$ to $n-1$ in the right had side of the formula \mythetag{5.11}.
Performing these calculations repeatedly, we derive the following formula:
$$
\hskip -2em
\int\limits^{+\infty}_y\!\!\varPhi_n(t)\,A_n(y/t)\,dt=-\!\!\int\limits^{+\infty}_y
\!\!\varphi{\kern 0.7 pt'}(t)\,\,dt.
\mytag{5.13}
$$
Due to \mythetag{5.1} the integral in the right hand side of \mythetag{5.13} is 
transformed to 
$$
\hskip -2em
-\!\!\int\limits^{+\infty}_y
\!\!\varphi{\kern 0.7 pt'}(t)\,\,dt=-\varphi(t)\,{\vrule height 16pt depth 
7pt}^{\,\infty}_{\ y}=\,\varphi(y).
\mytag{5.14}
$$
Combining the formulas \mythetag{5.13} and \mythetag{5.14} we derive the required 
formula \mythetag{5.4}. Thus, the lemma~\mythelemma{5.1} is proved for $n\geqslant 2$.
\par 
     If $n=1$ the formula \mythetag{5.12} coincides with \mythetag{5.13}. Similarly,
if $n=0$ the formula \mythetag{5.11} coincides with \mythetag{5.13}. Therefore, the
formula \mythetag{5.4} is proved for all $n\geqslant 0$. The proof of the 
lemma~\mythelemma{5.1} is over.\qed\enddemo
     Now, in addition to \mythetag{5.1} and \mythetag{5.2}, assume that the function
$\varphi(t)$ satisfies the following auxiliary condition at the point $t=0$\,:
$$
\lim_{t\to+0}\left(t^{s+\omega}\,\varphi^{(s)}\right)=0\text{\ \ for any \ }
\omega>0\text{\ \ and for all \ }s=0,\,\ldots,\,n+2.\quad
\mytag{5.15}
$$
Applying \mythetag{5.15} to \mythetag{5.3}, we find that $\varPhi_n(t)$ satisfies
the condition
$$
\hskip -2em
\lim_{t\to+0}\left(t^{1+\omega}\,\varPhi_n(t)\right)=0\text{\ \ for any \ }
\omega>0.
\mytag{5.16}
$$
\mylemma{5.2} For any $0<\alpha<\varepsilon$ the function \mythetag{5.3} satisfies 
the integral relationship
$$
\hskip -2em
\int\limits^{+\infty}_0\!\!\varPhi_n(t)\,\,t^\alpha\,dt=
-\alpha\prod^n_{k=1}\Bigl(1+\frac{\alpha}{k}\Bigr)
\!\int\limits^{+\infty}_0
\!\!t^\alpha\,\varphi{\kern 0.7 pt'}(t)\,dt
\text{\ \ \ for \ \ }n\geqslant 0,
\mytag{5.17}
$$
provided the conditions \mythetag{5.1}, \mythetag{5.2}, and \mythetag{5.15} are fulfilled. 
\endproclaim
     Assume that the conditions \mythetag{5.1}, \mythetag{5.2}, and \mythetag{5.15} are 
fulfilled. Then the function $\varPhi_n(t)$ satisfies the condition \mythetag{5.6}. If
$\alpha<\varepsilon$, the condition \mythetag{5.6} means that the integral in the left
hand side of the equality \mythetag{5.17} converges at infinity. The integral in the
right hand side of the equality \mythetag{5.17} also converges at infinity due to the
condition \mythetag{5.2} and the inequality $\alpha<\varepsilon$.\par
     The condition \mythetag{5.15} leads to the condition \mythetag{5.16} for the
function $\varPhi_n(t)$. Hence from $\alpha>0$ we derive that the integral in the left
hand side of the equality \mythetag{5.17} converges at the point $t=0$. Similarly the
integral in the right hand side of the equality \mythetag{5.17} converges at $t=0$ due
to \mythetag{5.15} and the inequality $\alpha>0$. 
\demo{Proof of the lemma~\mythelemma{5.2}} The proof is pure calculations. \pagebreak 
Upon substituting the formula \mythetag{5.3} into the left hand side of \mythetag{5.17} 
we can integrate by parts:
$$
\gather
\int\limits^{+\infty}_0\!\!\varPhi_n(t)\,\,t^\alpha\,dt=
-\!\!\int\limits^{+\infty}_0\!\!\frac{d}{dt}\left(\frac{(-t)^{n+1}}{n!}
\,\,\varphi^{(n+1)}(t)\right)t^\alpha\,dt=\\
=-\frac{(-1)^{n+1}}{n!}\,t^{\alpha+n+1}\,\varphi^{(n+1)}(t)\,{\vrule height 16pt 
depth 7pt}^{\,\infty}_{\,0}\!-\,\alpha\,\frac{(-1)^n}{n!}
\!\!\int\limits^{+\infty}_0\!\!t^{\alpha+n}\,\varphi^{(n+1)}(t)\,dt.
\endgather
$$
Note that $t^{\alpha+n+1}\,\varphi^{(n+1)}(t)\to 0$ as $t\to 0$ due to \mythetag{5.15}
since $\alpha>0$. Similarly $t^{\alpha+n+1}\,\varphi^{(n+1)}(t)\to 0$ as $t\to +\infty$ 
due to \mythetag{5.2} since $\alpha<\varepsilon$. As a result we get
$$
\hskip -2em
\int\limits^{+\infty}_0\!\!\varPhi_n(t)\,\,t^\alpha\,dt=-\alpha\,\frac{(-1)^n}{n!}
\!\!\int\limits^{+\infty}_0\!\!t^{\alpha+n}\,\varphi^{(n+1)}(t)\,dt.
\mytag{5.18}
$$\par
     In order to transform \mythetag{5.18} we continue integrating by parts and we shall 
decrease by $1$ the order of derivatives in each step. In the first step we get
$$
\gather
\int\limits^{+\infty}_0\!\!\varPhi_n(t)\,\,t^\alpha\,dt=-\alpha\,\frac{(-1)^n}{n!}
\,t^{\alpha+n}\,\varphi^{n}(t)\,{\vrule height 16pt 
depth 7pt}^{\,\infty}_{\,0}-\alpha\,\frac{(-1)^{n-1}}{n!}\,(\alpha+n)
\!\!\int\limits^{+\infty}_0\!\!t^{\alpha+n-1}\,\varphi^{(n)}(t)\,dt.
\endgather
$$
The non-integral term in the above formula vanishes due to \mythetag{5.2} and \mythetag{5.15}
since $\alpha>0$ and $\alpha<\varepsilon$. Hence the formula \mythetag{5.18} transforms to
$$
\hskip -2em
\int\limits^{+\infty}_0\!\!\varPhi_n(t)\,\,t^\alpha\,dt=-\alpha\,\frac{(-1)^{n-1}}{(n-1)!}
\left(1+\frac{\alpha}{n}\right)
\!\!\int\limits^{+\infty}_0\!\!t^{\alpha+n-1}\,\varphi^{(n)}(t)\,dt.
\mytag{5.19}
$$
Comparing the right hand sides of \mythetag{5.18} and \mythetag{5.19}, we see that they 
differ by passing from $n$ to $n-1$ and in \mythetag{5.19} we have the additional factor 
$(1-\alpha/n)$. Performing the above procedure repeatedly we shall reduce $n$ to $0$ and
gain more additional factors. They form the following product 
$$
\hskip -2em
\left(1+\frac{\alpha}{n}\right)\,\left(1+\frac{\alpha}{n-1}\right)\cdot\ldots\cdot
\left(1+\frac{\alpha}{1}\right)=\prod^n_{k=1}\left(1+\frac{\alpha}{k}\right).
\mytag{5.20}
$$
Due to \mythetag{5.20} the formula \mythetag{5.18} reduces to \mythetag{5.17}. Thus, we
have proved the lemma~\mythelemma{5.2} for $n\geqslant 2$.\par
     If $n=0$, the formula \mythetag{5.18} is equivalent to \mythetag{5.17}. Similarly,
if $n=1$ the formula \mythetag{5.19} is equivalent to \mythetag{5.17}. The 
lemma~\mythelemma{5.2} is proved.
\qed\enddemo
\head
6. Application to Khabibullin's conjecture.
\endhead
     The integral inequalities \mythetag{1.3} and \mythetag{1.4} in Khabibullin's
conjecture~\mytheconjecture{1.2} are related to some special choice of the function 
$\varphi(t)$ in \mythetag{5.4} and \mythetag{5.17}. Let's recall that the inequality 
\mythetag{1.3} now is written as \mythetag{3.2} in terms of the kernel function 
$A_{n-1}(x)$. Upon changing the variable of integration $x$ for $y=t\,x$ in 
\mythetag{3.2} this inequality transforms to the following one:
$$
\hskip -2em
\int\limits^{\,t}_0  A_{n-1}(y/t)\,q(y)\,\frac{dy}{t}\leqslant t^{\alpha-1}.
\mytag{6.1}
$$
Since $t\geqslant 0$ in \mythetag{6.1}, this inequality can be written as 
$$
\hskip -2em
\int\limits^{\,t}_0  A_{n-1}(y/t)\,q(y)\,dy\leqslant t^\alpha.
\mytag{6.2}
$$
Assume that $\varPhi_{n-1}(t)$ is positive on the interval $(0,+\infty)$ for some choice 
of $\varphi(t)$ in \mythetag{5.3}. Then we can multiply both sides of \mythetag{6.2} by
$\varPhi_{n-1}(t)$ and integrate over $t$ from $0$ to $+\infty$. As a result we derive the
following inequality:
$$
\hskip -2em
\int\limits^{+\infty}_0\!\!\varPhi_{n-1}(t)\!\left(\,\int\limits^{\,t}_0\!A_{n-1}(y/t)
\,q(y)\,dy\!\right)\!dt
\,\,\leqslant\!\int\limits^{+\infty}_0\!\!\varPhi_{n-1}(t)\,t^\alpha\,dt.
\mytag{6.3}
$$
Upon changing the order of integration in the left hand side of \mythetag{6.3} we get
$$
\hskip -2em
\int\limits^{+\infty}_0\left(\int\limits^{+\infty}_y\!\varPhi_{n-1}(t)\,A_{n-1}(y/t)
\,dt\!\right)\!q(y)\,dy
\,\,\leqslant\!\int\limits^{+\infty}_0\!\!\varPhi_{n-1}(t)\,t^\alpha\,dt.
\mytag{6.4}
$$
Now we can apply \mythetag{5.4} and \mythetag{5.17} to \mythetag{6.4}. This yields
$$
\hskip -2em
\int\limits^{+\infty}_0\!\!\varphi(y)\,q(y)\,dy
\,\,\leqslant-\alpha\prod^{n-1}_{k=1}\Bigl(1+\frac{\alpha}{k}\Bigr)
\!\int\limits^{+\infty}_0
\!\!t^\alpha\,\varphi{\kern 0.7 pt'}(t)\,dt
\mytag{6.5}
$$\par
     The left hand side of \mythetag{6.5} is an integral depending on $q(t)$, 
while its right hand side is a number. Comparing \mythetag{6.5} and \mythetag{1.4}, 
we see that our proper choice is
$$
\hskip -2em
\varphi(t)=\ln(1+t^{-2\,\alpha})\text{, \ where \ }\alpha>0.
\mytag{6.6}
$$
It is easy to verify that the function \mythetag{6.6} satisfies the condition 
\mythetag{5.15}. Moreover, it satisfies the conditions \mythetag{5.1} and
\mythetag{5.2} for $\varepsilon=2\,\alpha$.\par
     The integral in the right hand side of the inequality \mythetag{6.5} for the 
function \mythetag{6.6} can be calculated explicitly. Indeed, we have
$$
\pagebreak
\hskip -2em
\gathered
\int\limits^{+\infty}_0\!\!t^\alpha\,\varphi{\kern 0.7 pt'}(t)\,dt\,=
\!\!\int\limits^{+\infty}_0\!\!t^\alpha\,\frac{-2\,\alpha\,t^{-2\,\alpha-1}}
{1+t^{-2\,\alpha}}\,dt=-2\!\!\int\limits^{+\infty}_0\!\!\frac{\alpha\,t^{\alpha-1}}
{1+t^{2\,\alpha}}\,dt=\\
=-2\!\!\int\limits^{+\infty}_0\!\!\frac{(t^\alpha)'\,dt}{1+(t^\alpha)^2}
=-2\!\!\int\limits^{+\infty}_0\!\!\frac{dz}{1+z^2}=-2\,\arctan(z)\,{\vrule height 16pt 
depth 7pt}^{\,\infty}_{\ 0}=-\pi.
\endgathered
\mytag{6.7}
$$
Applying \mythetag{6.6} and \mythetag{6.7} to \mythetag{6.5}, we get the inequality 
coinciding with the required inequality \mythetag{1.4} in Khabibullin's conjecture.
\head
7. The analysis of the transition function.
\endhead
     The result of the previous section shows that Khabibullin's 
conjecture~\mytheconjecture{1.2} is valid, provided $\varPhi_{n-1}(t)\geqslant 0$ for all 
$t\in(0,+\infty)$. Unfortunately, the transition function $\varPhi_{n-1}(t)$ produced by
the function \mythetag{6.6} is not always positive.\par
     Note that the function \mythetag{6.6} depends on the parameter $\alpha>0$ coinciding
with the parameter $\alpha$ in Khabibullin's conjecture~\mytheconjecture{1.2}. Therefore
the transition function $\varPhi_{n-1}(t)$ depends on two parameters $n$ and $\alpha$, i\.
\,e\. $\varPhi_{n-1}(t)=\varPhi_{n-1}(\alpha,t)$. Our goal in this section is to search
some values of $n$ and $\alpha$ for which 
$$
\hskip -2em
\varPhi_{n-1}(\alpha,t)\geqslant 0\text{\ \ for all \ }t>0.
\mytag{7.1}
$$\par
     Let's begin with $n=1$. Then $\varPhi_{n-1}(\alpha,t)=\varPhi_0(\alpha,t)$. From
\mythetag{5.3} and \mythetag{6.6} for this case we derive the following expression for 
$\varPhi_0(\alpha,t)$:
$$
\hskip -2em
\gathered
\varPhi_0(\alpha,t)=\left(t\,\,\varphi{\kern 0.7 pt'}(t)\right)'
=\frac{d}{dt}\!\left(\frac{-2\,\alpha\,t^{-2\,\alpha}}
{1+t^{-2\,\alpha}}\right)=\\
=\frac{d}{dt}\!\left(\frac{-2\,\alpha}{1+t^{2\,\alpha}}\right)
=\frac{4\,\alpha^2}{t}\cdot\frac{t^{2\,\alpha}}{(1+t^{2\,\alpha})^2}.
\endgathered
\mytag{7.2}
$$
It is easy to see that the function \mythetag{7.2} is positive for all $t>0$, i\.\,e\.
the condition \mythetag{7.1} is fulfilled. As a result we get the following theorem.
\mytheorem{7.1} If\/ $n=1$, Khabibullin's conjecture~\mytheconjecture{1.2} is valid 
for all $\alpha>0$.
\endproclaim
     Let's proceed to the case $n=2$. In this case $\varPhi_{n-1}(\alpha,t)
=\varPhi_1(\alpha,t)$. Substituting the function \mythetag{6.6} into the formula 
\mythetag{5.3}, we derive
$$
\hskip -2em
\varPhi_1(\alpha,t)=-\left(t^2\,\varphi{\kern 0.7 pt''}(t)\right)'
=\frac{d}{dt}\!\left(\!t^2\,\frac{d}{dt}\!\left(\frac{2\,\alpha\,t^{-2\,\alpha-1}}
{1+t^{-2\,\alpha}}\right)\!\!\right)\!.
\mytag{7.3}
$$
Upon expanding \mythetag{7.3} and simplifying the obtained expression, we get
$$
\hskip -2em
\varPhi_1(\alpha,t)=\frac{4\,\alpha^2}{t}\cdot
\frac{(2\,\alpha+1)\,t^{4\,\alpha}+(1-2\,\alpha)
\,t^{2\,\alpha}}{(1+t^{2\,\alpha})^3}.
\mytag{7.4}
$$
The formula \mythetag{7.4} can be written as 
$$
\hskip -2em
\varPhi_1(\alpha,t)=\frac{4\,\alpha^2}{t}\cdot
\frac{t^{2\,\alpha}\,P_1(\alpha,t^{2\,\alpha})}{(1+t^{2\,\alpha})^3},
\mytag{7.5}
$$
where $P_1=P_1(\alpha,z)$ is a polynomial of the variable $z=t^{2\,\alpha}$:
$$
\hskip -2em
P_1(\alpha,z)=(2\,\alpha+1)\,z+(1-2\,\alpha).
\mytag{7.6}
$$
The  polynomial $P_1(z)$ in \mythetag{7.6} is positive for all $z>0$ if and only 
if \,$0<\alpha\leqslant 1/2$. Hence we have the following theorem. 
\mytheorem{7.2} If\/ $n=2$, Khabibullin's conjecture~\mytheconjecture{1.2} is valid 
for all\/ $0<\alpha\leqslant 1/2$.
\endproclaim
     The next case is $n=3$. In this case $\varPhi_{n-1}(\alpha,t)
=\varPhi_2(\alpha,t)$. Substituting the function \mythetag{6.6} into the formula 
\mythetag{5.3}, we derive
$$
\hskip -2em
\varPhi_2(\alpha,t)=\frac{\left(t^3\,\varphi{\kern 0.7 pt'''}(t)\right)'}{2}
=-\frac{d}{dt}\!\left(\!t^3\,\frac{d^{\kern 0.4pt 2}}{dt^2}\!\left(\frac{\alpha
\,t^{-2\,\alpha-1}}{1+t^{-2\,\alpha}}\right)\!\!\right)\!.
\mytag{7.7}
$$
It is preferable to write the formula \mythetag{7.7} as follows:
$$
\hskip -2em
\varPhi_2(\alpha,t)=-\frac{d}{dt}\!\left(\!t^3\,\frac{d^{\kern 0.4pt 2}}{dt^2}
\!\left(\frac{\alpha\,t^{-1}}{1+t^{2\,\alpha}}\right)\!\!\right)\!.
\mytag{7.8}
$$
By mens of direct calculations we transform the formula \mythetag{7.8} to 
$$
\hskip -2em
\varPhi_2(\alpha,t)=\frac{4\,\alpha^2}{t}\cdot
\frac{t^{2\,\alpha}\,P_2(\alpha,t^{2\,\alpha})}{(1+t^{2\,\alpha})^4},
\mytag{7.9}
$$
where $P_2=P_2(\alpha,z)$ is the following quadratic polynomial of the variable 
$z=t^{2\,\alpha}$:
$$
P_2(\alpha,z)=(\alpha+1)\,(1+2\,\alpha)\,z^2-2\,(2\,\alpha-1)\,(1+2\,\alpha)\,z
+(2\,\alpha-1)(\alpha-1).\quad
\mytag{7.10}
$$
\mylemma{7.1} Let $P(z)=A\,z^2+B\,z+C$ be a general quadratic polynomial with real
coefficients $A$, $B$, and $C$. Then $P(z)\geqslant 0$ for all $z>0$ if and only 
if one of the following three conditions is fulfilled:
\roster
\rosteritemwd=5pt
\item"1)" $A>0$, \ $B<0$, \ $B^2-4\,A\,C\leqslant 0$;
\item"2)" $A>0$, \ $B\geqslant 0$, \ $C\geqslant 0$;
\item"3)" $A=0$, \ $B\geqslant 0$, \ $C\geqslant 0$.
\endroster
\endproclaim
\demo{Proof} If $A\neq 0$, then $P(z)\to -\infty$ as $z\to +\infty$ for $A<0$ and 
$P(z)\to +\infty$ as $z\to +\infty$ for $A>0$. Therefore $A\geqslant 0$ is a necessary 
condition for $P(z)\geqslant 0$ for all $z>0$. If $A>0$, the graph of the function $P(z)$ 
is a parabola (see Fig\.~7.1 and Fig\.~7.2). If $A=0$, it is a straight line 
(see Fig\.~7.3). \vadjust{\vskip 45pt\hbox to 0pt{\kern 5pt \includegraphics{khab02.eps}\hss}\vskip 120pt}In the last case $P(z)\geqslant 0$ for all $z>0$ if and 
only if the following two inequalities are fulfilled:
$$
\xalignat 2
&\hskip -2em
C\geqslant 0,
&&B\geqslant 0.
\mytag{7.11}
\endxalignat
$$\par
      In the case of a parabolic graph, i\.\,e\. if $A>0$, the function $P(z)$ decreases 
for $z<z_{\sssize\text{min}}$ and $P(z)$ increases for $z>z_{\sssize\text{min}}$. Hence
there are two subcases where $P(z)\geqslant 0$ for all $z>0$. The first subcase is given 
by the inequalities (see Fig\.~7.1):
$$
\xalignat 2
&\hskip -2em
z_{\sssize\text{min}}>0,
&&P_{\sssize\text{min}}\geqslant 0.
\mytag{7.12}
\endxalignat 
$$
The second subcase corresponds to Fig\.~7.2. It is given by the inequalities
$$
\xalignat 2
&\hskip -2em
z_{\sssize\text{min}}\leqslant 0,
&&P(0)\geqslant 0.
\mytag{7.13}
\endxalignat 
$$
The elementary calculus yields
$$
\xalignat 3
&\hskip -2em
z_{\sssize\text{min}}=-\frac{B}{2\,A},
&&P_{\sssize\text{min}}=C-\frac{B^2}{4\,A},
&&P(0)=C.
\mytag{7.14}
\endxalignat 
$$
Applying \mythetag{7.14} to \mythetag{7.12} and taking into account that $A>0$,
we obtain
$$
\xalignat 2
&\hskip -2em
B<0,
&&B^2-4\,A\,C\leqslant 0.
\mytag{7.15}
\endxalignat 
$$
Similarly, applying \mythetag{7.14} to \mythetag{7.13} and taking into account that 
$A>0$, we get 
$$
\xalignat 2
&\hskip -2em
B\geqslant 0,
&&C\geqslant 0.
\mytag{7.16}
\endxalignat 
$$
Now in order to complete the proof of the lemma~\mythelemma{7.1} it is sufficient to 
complement the inequalities \mythetag{7.15} and \mythetag{7.16} with the inequality 
$A>0$ and to complement the inequalities \mythetag{7.11} with the equality $A=0$. 
\qed\enddemo
     having proved the lemma~\mythelemma{7.1}, we apply it to the polynomial
\mythetag{7.10}. In this case $A=(\alpha+1)\,(1+2\,\alpha)$. Since $\alpha>0$ in the
conjecture~\mytheconjecture{1.2}, we have $A>0$, i\.\,e\. the third option of the
lemma~\mythelemma{7.1} can be dropped from our further considerations. Since
$B=-2\,(2\,\alpha-1)\,(1+2\,\alpha)$ and $C=(2\,\alpha-1)(\alpha-1)$, the second 
option of the lemma~\mythelemma{7.1} combined with $\alpha>0$, leads to the following 
system of inequalities:
$$
\hskip -2em
\left\{\aligned
&(1-2\,\alpha)\,(1+2\,\alpha)\geqslant 0,\\
&(2\,\alpha-1)(\alpha-1)\geqslant 0,\\
&\alpha>0.\endaligned\right.
\mytag{7.17}
$$ 
The inequalities \mythetag{7.17} resolve to $0<\alpha\leqslant 1/2$.\par
     Let's proceed to the first option of the lemma~\mythelemma{7.1}. By means of direct 
calculations we get $B^2-4\,A\,C=12\,\alpha^2\,(2\,\alpha-1)(1+2\,\alpha)$. Therefore
the first option of the lemma~\mythelemma{7.1} combined with the inequality $\alpha>0$ 
yields
$$
\hskip -2em
\left\{\aligned
&(1-2\,\alpha)\,(1+2\,\alpha)<0,\\
&\alpha^2\,(2\,\alpha-1)(1+2\,\alpha)\leqslant 0,\\
&\alpha>0.\endaligned\right.
\mytag{7.18}
$$ 
It is easy to find that the inequalities \mythetag{7.18} are mutually contradictory. They 
cannot be satisfied simultaneously. As a result we conclude that the values of the polynomial 
\mythetag{7.10} are non-negative for all $z>0$ \pagebreak if and only 
if\/ $0<\alpha\leqslant 1/2$. Applying this fact to the transition function \mythetag{7.9}, 
we get the following theorem. 
\mytheorem{7.3} If\/ $n=3$, Khabibullin's conjecture~\mytheconjecture{1.2} is valid 
for all\/ $0<\alpha\leqslant 1/2$.
\endproclaim
\head
8. Recurrent formulas for the transition functions. 
\endhead
     The transition function $\varPhi_n(\alpha,t)$ introduced in Section~5 is given by
the explicit formula \mythetag{5.3}. However, for our purposes we need to have the 
recurrent formula expressing $\varPhi_n(\alpha,t)$ through the function 
$\varPhi_{n-1}(\alpha,t)$. Here is this formula
$$
\hskip -2em
\varPhi_n=-\frac{t^n}{n}\,\frac{d}{dt}\!\left(\frac{\varPhi_{n-1}}{t^{n-1}}\!\right).
\mytag{8.1}
$$
The formula \mythetag{8.1} is easily proved by means of direct calculations with the use 
of the initial formula \mythetag{5.3}.\par
      The transition functions $\varPhi_n(\alpha,t)$ associated with Khabibullin's conjecture~\mytheconjecture{1.2} correspond to the special choice \mythetag{6.6} of the 
function $\varphi(t)$ in \mythetag{5.3}. The formula \mythetag{7.2}, \mythetag{7.5}, 
and \mythetag{7.9} were derived exactly for this special choice of $\varphi(t)$. 
Comparing these three formulas with each other, we write 
$$
\hskip -2em
\varPhi_n(\alpha,t)=\frac{4\,\alpha^2}{t}\cdot
\frac{t^{2\,\alpha}\,P_n(\alpha,z)}{(1+t^{2\,\alpha})^{n+2}},
\mytag{8.2}
$$
where $z=t^{2\,\alpha}$ and $P_n(\alpha,z)$ is a polynomial of the degree $n$ with 
respect to the variable $z$. The formula \mythetag{8.2} is proved by induction on $n$. 
Using \mythetag{8.1} one can derive the following recurrent formula for the polynomials 
$P_n(z)$ in \mythetag{8.2}:
$$
\hskip -2em
\gathered
P_n(\alpha,z)=\left(\!(2\,\alpha+1)z+\left(1-\frac{2\,\alpha}{n}\right)\!\!\right)
P_{n-1}(\alpha,z)\,-\\
\vspace{1ex}
-\,\frac{2\,\alpha\,z\,(z+1)}{n}\,P^{\,\prime}_{n-1}(\alpha,z).
\endgathered
\mytag{8.3}
$$
The formula \mythetag{7.2} yields the base for applying the recurrent formula
\mythetag{8.3}:
$$
\hskip -2em
P_0(\alpha,z)=1.
\mytag{8.4}
$$
Combining the formulas \mythetag{8.3} and \mythetag{8.4}, one can calculate the 
polynomial $P_n(\alpha,z)$ explicitly for each particular $n$.
\mytheorem{8.1} If\/ $0<\alpha\leqslant 1/2$, then $P_n(\alpha,z)\geqslant 0$ for all $z>0$.
\endproclaim
\mytheorem{8.2} If\/ $0<\alpha\leqslant 1/2$, then $\varPhi_n(\alpha,t)\geqslant 0$ for 
all $t>0$.
\endproclaim
     Due to \mythetag{8.2} the theorems~\mythetheorem{8.1} and \mythetheorem{8.2} are
equivalent to each other. However, the theorem~\mythetheorem{8.2} is easier to prove.
\demo{Proof of the theorem \mythetheorem{8.2}} Applying the recurrent formula
\mythetag{8.1} we immediately derive the following expression for $\varPhi_n(\alpha,t)$:
$$
\hskip -2em
\varPhi_n(\alpha,t)=\frac{(-1)^n}{n!}\,\frac{d^{\kern 0.3pt n}\varPhi_0(\alpha,t)}{dt^n}.
\mytag{8.5}
$$
The function $\varPhi_0(\alpha,t)$ is given by the formula \mythetag{7.2}. We write it as
$$
\hskip -2em
\varPhi_0(\alpha,t)=\frac{4\,\alpha^2}{(1+t^{2\,\alpha})^2\,t^{1-2\,\alpha}}.
\mytag{8.6}
$$
The function \mythetag{8.6} belongs to the class of functions $K\!K(\beta)$ with 
$\beta=2\,\alpha$. By definition this class of functions consists of all linear combinations 
of functions
$$
\hskip -2em
\psi(t)=\frac{1}{t^{\kern 0.2pt\gamma}\,(1+t^\beta)^k},\quad\gamma\geqslant 0, 
\quad k\geqslant 0, 
\mytag{8.7}
$$
with non-negative coefficients. The class of functions $K\!K(\beta)$ is similar to the class
of functions $K(\beta)$ defined in \mycite{4}.\par
      Obviously, each function $f\in K\!K(\beta)$ is non-negative, i\.\,e\. $f(t)\geqslant 0$ 
for all $t>0$. For $0<\beta\leqslant 1$ the class of functions $K\!K(\beta)$ is closed with 
respect to the operator
$$
\hskip -2em
D=-\frac{d}{dt}. 
\mytag{8.8}
$$
Indeed, applying the operator \mythetag{8.8} to the function \mythetag{8.7}, we get
$$
\hskip -2em
D\psi(t)=\gamma\cdot\frac{1}{t^{\kern 0.2pt\gamma+1}\,(1+t^\beta)^k}
+(k\,\beta)\cdot\frac{1}{t^{\kern 0.2pt\gamma+1-\beta}\,(1+t^\beta)^{k+1}}
\mytag{8.9}
$$
As we see, for $0<\beta\leqslant 1$ the right hand side of \mythetag{8.9} is a linear 
combination of two functions of the form \mythetag{8.7} with two non-negative coefficients 
$\gamma$ and $k\,\beta$.\par
     Now let's return back to the formulas \mythetag{8.5} and \mythetag{8.6}. In terms
of the operator \mythetag{8.8} the formula \mythetag{8.5} is written as follows:
$$
\hskip -2em
\varPhi_n(\alpha,t)=\frac{D^n\varPhi_0(\alpha,t)}{n!}.
\mytag{8.10}
$$
Since $\beta=2\,\alpha$ and $0<\alpha\leqslant 1/2$ is equivalent to $0<\beta\leqslant 1$,
from $\varPhi_0(\alpha,t)\in K\!K(\beta)$ and \mythetag{8.10} we derive $\varPhi_n(\alpha,t)
\in K\!K(\beta)$. Hence $\varPhi_n(\alpha,t)\geqslant 0$ for all $t>0$. This means that the 
theorem~\mythetheorem{8.2} is proved. \qed\enddemo
     From the theorem~\mythetheorem{8.2} we immediately derive the following theorem.
\mytheorem{8.3} Khabibullin's conjecture~\mytheconjecture{1.2} is valid 
for all\/ $0<\alpha\leqslant 1/2$ and for all integer $n>0$.
\endproclaim
\head
9. Conclusions.
\endhead
     The theorem~\mythetheorem{8.3} is not a new result. It is known from \mycite{4}.
However, using the transition function $\varPhi_{n-1}(\alpha,t)$, we define the integral 
transformation 
$$
\hskip -2em
\psi(t)\mapsto\int\limits^{+\infty}_0\!\!\varPhi_{n-1}(t)\,\psi(t)\,dt.
\mytag{9.1}
$$
Due to the lemma~\mythelemma{5.1} the transformation \mythetag{9.1} can be worth in 
studying Khabi\-bullin's conjecture~\mytheconjecture{1.2} for its parameters $n$ and 
$\alpha$ in ranges beyond those covered by the theorem~\mythetheorem{8.3}. 


\Refs
\ref\myrefno{1}\by Khabibullin~B.~N.\paper Paley problem for plurisubharmonic functions 
of a finite lower order\jour Mat\. Sbornik\vol 190\issue 2\yr 1999\pages 145-157
\endref
\ref\myrefno{2}\by Khabibullin~B.~N.\paper The representation of a meromorphic function as 
a quotient of entire functions and the Paley problem in $\Bbb C^n$: survey of some results
\jour Mathematical Physics, Analysis, and geometry (Ukraine) \yr 2002\vol 9\issue 2
\pages 146-167\moreref see also
\myhref{http://arxiv.org/abs/math.CV/0502433}{math.CV/0502433} in Electronic Archive 
\myEarXivlink
\endref
\ref\myrefno{3}\by Khabibullin~B.~N.\paper A conjecture on some estimates for integrals
\publ e-print \myhref{http://arXiv.org/abs/1005.3913/}{arXiv:1005.3913} in Electronic 
Archive \myEarXivlink
\endref
\ref\myrefno{4}\by Baladai~R.~A, Khabibullin~B.~N.\paper Three equivalent conjectures on 
an estimate of integrals\publ e-print \myhref{http://arXiv.org/abs/1006.5140/}{arXiv:1006.5140} 
in Electronic Archive \myEarXivlink
\endref
\ref\myrefno{5}\by Zorich~V.~A.\book Mathematical analysis. 
\rm Part~\uppercase\expandafter{\romannumeral 2}\publ Nauka publishers\publaddr
Moscow\yr 1984
\endref
\endRefs
\enddocument
\end